\documentclass{article}

\usepackage{amsmath,amssymb,amsthm}
\usepackage{cite}
\usepackage{a4wide}
\usepackage{color}
\usepackage{txfonts}
\usepackage{subfigure, verbatim}
\usepackage{mathptmx}
\usepackage{rotating}
\usepackage{graphicx}
\usepackage{booktabs}
\usepackage{float}
\DeclareMathAlphabet{\mathpzc}{OT1}{pzc}{m}{it}

\newtheorem{theorem}{Theorem}[section]
\newtheorem{lemma}[theorem]{Lemma}
\newtheorem{corollary}[theorem]{Corollary}

\theoremstyle{definition}
\newtheorem{definition}[theorem]{Definition}

\begin{document}
\author{Samia Bashir$^1$, Amjad Hussain$^1$}

\title{\bf{Bounds for $p$-adic Hardy-type Operators and Commutator On $p$-adic Variable Herz-Morrey Spaces }}%Variable Herz-Morrey spaces with Hardy-type operators

\date{}
\maketitle

% --------------------------------------------

\begin{center}
$^{1}$ Department of Mathematics, Quaid-I-Azam University 45320, Islamabad 44000, Pakistan\\
%$^{2}$ Department of Mathematics, College of Science Al-Zulfi, Majmaah University, Al-Majmaah, P.O Box 66, Majmaah 11952, Saudi Arabia.\\
%$^{3}$ Department of Mathematics, Bonga University, Bonga, Ethiopia
%
%$^*$Corresponding Author email: ilyaskhan@tdtu.edu.vn
email: sbashir@math.qau.edu.pk
\end{center}

% --------------------------------------------

\begin{abstract}
This paper showed that fractional p-adic Hardy operator norms in p-adic Herz-Morrey spaces with varying exponents are bounded. Corresponding commutator operators are also estimated for p-adic variable central bounded mean oscillations (CBMO).
%In this paper, we showed that fractional p-adic Hardy operators norms in p-adic Herz-Morrey spaces with varying exponents are bounded. Additionally corresponding commutator operators are estimated for p-adic variable central bounded mean oscillations (CBMO).
% This article establishes the boundedness of variable exponent p-adic Herz spaces and variable exponent p-adic Morrey-Herz spaces for the fractional p-adic Hardy operator and its adjoint operator. On the aforementioned spaces, we also obtain the bounds for the commutators of fractional p-adic Hardy operators.
\bigskip


\bigskip

\end{abstract}

% ---------------------------------------------
\section{Introduction}According to the well-known Ostrowski theorem \cite{C5.1}, any nontrivial valuation on the field of rational numbers $\mathbb{Q}$ is equivalent either to the p-adic valuation $|\cdot|_{p}$ or to one of real valuations $|\cdot|$, where p is a prime number. The former norm is defined as follows: if any rational number $x\neq0$ is denoted as $x=p^{\gamma}s/t$, where $\gamma=\gamma(x)\in\mathbb{Z}$ and the integers s, t are not divisible by p, then $|x|_{p}=p^{\gamma}$, $|0|_{p}=0$. The norm $|\cdot|_{p}$ satisfies the strong triangle inequality $|x+y|_{p}\leq max(|x|_{p},|y|_{p})$. The extended form of any $x\neq0\in \mathbb{Q}_{p}$ (field of p-adic numbers) is given in \cite{C5.2} as
\begin{equation}\label{e5.1}
x= p^{\gamma}\sum_{i=0}^{\infty}x_{i}p^{i},
\end{equation}
where $x_{i},\gamma\in\mathbb{Z}$,\hspace{0.2cm}$x_{i}\in\frac{\mathbb{Z}}{p\mathbb{Z}_{p}},\hspace{0.2cm} x_{0}\neq0.$\\

In what follows the n-dimensional vector space $\mathbb{Q}_{p}^{n}$ over $\mathbb{Q}_{p}$ is equipped with the following absolute value:
\begin{equation}\label{e5.2}
|x|_{p}=\max\limits_{1\leq j\leq n}|x_{j}|_{p}.
\end{equation}
Let $B_{\gamma}(\mathfrak{a})$ and $S_{\gamma}(\mathfrak{a})$ represent respectively the ball and sphere of $\mathbb{Q}_{p}^{n}$ centered at $\mathfrak{a}\in \mathbb{Q}_{p}^{n}$ and radius $p^{\gamma}>0$:
\begin{equation}\begin{aligned}[b]\label{e5.3}
B_{\gamma}(\mathfrak{a})=\{\mathfrak{x}\in \mathbb{Q}_{p}^{n}:|\mathfrak{x}-\mathfrak{a}|_{p}\leq p^{\gamma}\},
\end{aligned}\end{equation}
\begin{equation}\begin{aligned}[b]\label{e5.4}
S_{\gamma}(\mathfrak{a})=\{\mathfrak{x}\in \mathbb{Q}_{p}^{n}:|\mathfrak{x}-\mathfrak{a}|_{p}= p^{\gamma}\}=B_{\gamma}(\mathfrak{a})\setminus B_{\gamma-1}(\mathfrak{a}).
\end{aligned}\end{equation}
We denote that $B_{\gamma}(0)=B_{\gamma}$ and $S_{\gamma}(0)=S_{\gamma}$. Since $\mathbb{Q}_{p}^{n}$ is a locally compact commutative group with respect to addition, in $\mathbb{Q}_{p}^{n}$ there exists a positive Haar measure $d\mathfrak{x}$ under shift $d(\mathfrak{x}+\mathfrak{a})=d\mathfrak{x}$. It is worth noting that $d\mathfrak{x}$ is normalized by the equality $\int_{B_{0}}d\mathfrak{x}=1$. It is easy to find $\int_{B_{\gamma}(\mathfrak{a})}d\mathfrak{x}=p^{n\gamma}$ and $\int_{S_{\gamma}(\mathfrak{a})}d\mathfrak{x}=p^{n\gamma}(1-p^{-n})$ for any $\mathfrak{a}\in \mathbb{Q}_{p}^{n}$.\\
The field of p-adic numbers can be applied in many scientific fields. In physics, the groundbreaking application is the p-adic AdS/CFT \cite{C5.3}. Khrennikov et al. \cite{C5.4} enhanced p-adic wavelet for modeling reaction-diffusion dynamics. Its application in biology includes the models for hierarchical structures of genetic code \cite{C5.5} and protein \cite{C5.6}. Furthermore, p-adic numbers have found a novel application in harmonic analysis and mathematical physics (see, for example, \cite{C5.2,C5.7,C5.8,C5.9,C5.10,C5.11} and related references).\\
The topic of fractional calculus is undergoing fast development with more and more appealing applications in the real world (see, for instance, \cite{C5.12, C5.13, C5.14, C5.15, C5.16}). fractional integral operators are an important part of the mathematical analysis as they construct and formulate inequalities which have multiple applications in scientific areas that can be seen in the existing literature \cite{C5.17, C5.18, C5.19, C5.20}. In this sense, Wu \cite{C5.21} defined the p-adic fractional Hardy operator as:
\begin{equation}\label{e5.5}
H^{p}_{\beta}f(\mathfrak{x})=\frac{1}{|\mathfrak{x}|^{n-\beta}_{p}}\int_{|t|_{p}\leq|\mathfrak{x}|_{p}}f(\mathrm{t})dt,\hspace{0.2cm}\mathfrak{x}\in\mathbb{Q}_{p}^{n}\setminus\{0\},
\end{equation}
where $f\in L_{loc}(\mathbb{Q}_{p}^{n})$ and $0\leq\beta\leq$. Also, he gave the following definition of its commutators:
\begin{equation}\label{e5.6}
H^{p}_{\beta,b}f= bH^{p}_{\beta}f- H^{p}_{\beta}(bf).
\end{equation}
If $\beta=0$, the fractional p-adic Hardy type operator is the p-adic Hardy operator \cite{C5.22, C5.23}. Recently, the commutators of the Hardy operator, Hardy-Ces\`{a}ro operator, and Hausdorff operator have been extensively studied on the real field, p-adic field and Heisenberg group (see e.g., \cite{C5.24, C5.25, C5.26, C5.27, C5.28} and references therein for more details). As is well known, the theory of function spaces with variable exponents has some essential applications in the electronic fluid mechanics, recovery of graphics, elasticity, harmonic analysis, and partial differential equations (see e.g., [1, 2, 4, 5, 9, 12, 13, 19, 25, 26, 32–34, 36] and the references therein).

The Morrey spaces first appeared in 1938 in the work of Morrey [1] in relation to some problems in partial differential equations. In [2] the authors introduced central Morrey spaces. The Herz spaces are a class of function spaces introduced by Herz in the study of absolutely convergent Fourier transforms in 1968; see [3]. The complete theory of Herz spaces for the case of general indexes was established by Lu et al. in 2008; see [4]. Lu and Xu defined the homogeneous Morrey-Herz spaces in [5].

The theory of function spaces with variable exponent was extensively studied by researchers since the work of Kova´cik and R ˇ akosn ´ ´ık [6] appeared in 1991; see [7, 8] and the references therein. Many applications of these spaces were given, for example, in the modeling of electrorheological fluids [9], in the study of image processing [10], and in differential equations with nonstandard growth [11]. In 2009, Izuki established the Herz spaces with variable exponent and Morrey-Herz spaces with variable exponent; see [12, 13]. In [14], the authors introduced the nonhomogeneous central Morrey spaces of variable exponent. %Recently, Chacón-Cortés and Rafero gave new dimension to p-adic numbers by introducing p-adic variable exponent Lebesgue spaces in []. However, there is no theory of p-adic Herz-Morrey spaces with variable exponents in the field of p-adic numbers.

Chacón-Cortés and Rafero recently gave p-adic numbers a new dimension by introducing p-adic variable exponent Lebesgue spaces in []. However, in the study of p-adic numbers, a theory of p-adic Herz-Morrey spaces with variable exponents is missing. To bridge this gap, we introduce p-adic Herz-Morrey spaces with variable exponents and prove that p-adic fractional Hardy operators on these spaces are bounded.

The following is a summary of the paper. Section 2 provides a brief introduction to the study of p-adic function spaces with variable exponents and a review of several well-established results in this area. In Section 3, we look into the boundedness of p-adic fractional Hardy operators and their corresponding commutators in the context of p-adic variable Herz space, when the symbol functions are members of p-adic CBMO spaces with varying exponents. In Section 4, we talk about how the same boundedness holds for the p-adic fractional Hardy operator and its commutators on p-adic variable Morrey-Herz type spaces.

\section{Preliminaries}
In this section, we fix the notation and gather some fundamental findings on p-adic analysis that we will utilize throughout the paper.

\subsection{Some Function Spaces}

For any $x\in\mathbb{Q}^{n}_{p}$, a complex-valued function $\hbar$ defined on $\mathbb{Q}^{n}_{p}$ is referred to as a local constant if an integer $\ell(x)\in\mathbb{Z}$ exists such that
\begin{equation}\label{C2e10}
\hbar(x+x')=\hbar(x)\hspace{0.5cm}for\hspace{0.2cm}x'\in B_{\ell(x)}.
\end{equation}
If compactly supported $\hbar:\mathbb{Q}^{n}_{p}\longrightarrow\mathbb{C}$ is locally constant, it is called a Schwartz-Bruhat function (or a test function). $S(\mathbb{Q}^{n}_{p})=:S$ denotes the $\mathbb{C}$-vector space of Schwartz-Bruhat functions.\\
A measurable function $f:\mathbb{Q}^{n}_{p}\longrightarrow\mathbb{C}$ is a member of the Lebesgue space $L^{u}(\mathbb{Q}^{n}_{p})$, $1\leq u<\infty$, when
\begin{equation}\label{C2e11}
\|f\|^{u}_{L^{u}(\mathbb{Q}^{n}_{p})}=:\int_{\mathbb{Q}^{n}_{p}}|f(x)|^{u}dx<\infty,
\end{equation}
where
\begin{equation}\label{C2e12}
\int_{\mathbb{Q}^{n}_{p}}|f(x)|^{u}dx=:\lim\limits_{\gamma\rightarrow\infty}\int_{B_{\gamma}(0)}|f(x)|^{u}dx,
\end{equation}
if the limit does exist.\\
The concept of p-adic spaces of Lebesgue with a variable exponent is introduced here and also certain properties that will be needed in the next section are given; the proofs can be found in \cite{C2.19}.\\
If u is a function from $\mathbb{Q}^{n}_{p}$ to $[1,\infty)$, we call a variable exponent is a function that can be measured. The set of all variable exponents satisfying $u^{+}$ less than infinity is denoted by $\mathfrak{G}(\mathbb{Q}^{n}_{p})$, where $u^{+}=:esssup\{u(x):x\in\mathbb{Q}^{n}_{p}\}$ and $u^{-}=: essinf\{u(x):x\in\mathbb{Q}^{n}_{p}\}$.\\
We refer to this as the space of measurable functions $f:\mathbb{Q}^{n}_{p}\longrightarrow\mathbb{R}$ by $L^{u(\cdot)}(\mathbb{Q}^{n}_{p})$ for $u\in\mathfrak{G}(\mathbb{Q}^{n}_{p})$ such that
\begin{equation}\label{C2e13}
\|f\|_{L^{u(\cdot)}(\mathbb{Q}^{n}_{p})}=:inf\left\{\lambda>0:\wp_{u(\cdot)}\left(\frac{f}{\lambda}\right)\leq1\right\}<\infty,
\end{equation}
where $\wp_{u(\cdot)}(f)=:\int_{\mathbb{Q}^{n}_{p}}|f(x)|^{u(x)}dx$.

We now have the following for the Lebesgue space with variable exponent
\begin{equation}\label{C2e14}
\|f\|_{L^{u(\cdot)}(\mathbb{Q}^{n}_{p})}\leq\wp_{u(\cdot)}(f)+1,
\end{equation}
\begin{equation}\label{C2e15}
\wp_{u(\cdot)}(f)\leq\left(1+\|f\|_{L^{u(\cdot)}(\mathbb{Q}^{n}_{p})}\right)^{u^{+}},
\end{equation}
\begin{equation}\label{C2e16}
\|f\|_{L^{u(\cdot)}(\mathbb{Q}^{n}_{p})}=\||f|^{s}\|^{\frac{1}{s}}_{L^{u(\cdot)/s}(\mathbb{Q}^{n}_{p})},\hspace{0.5cm}s\in(0,u^{-}].
\end{equation}
In Lebesgue spaces where the exponent can vary, the Holder inequaity is true upto a multiplicative constant, i.e.
\begin{equation}\label{C2e17}
\int_{\mathbb{Q}^{n}_{p}}|f(x)g(x)|dx\leq C\|f\|_{L^{u(\cdot)}(\mathbb{Q}^{n}_{p})}\|g\|_{L^{u'(\cdot)}(\mathbb{Q}^{n}_{p})},
\end{equation}
here, $u$ and $u'$ are both conjugate exponents, which means that $1=1/u(x)+1/u'(x)$.\\
When there is a constant C that is positive for $u\in\mathfrak{G}(\mathbb{Q}^{n}_{p})$, we say that $u\in W_{0}(\mathbb{Q}^{n}_{p})$, as a result of which
\begin{equation}\label{C2e18}
\gamma\left(u^{-}(B_{\gamma}(x))-u^{+}(B_{\gamma}(x))\right)\leq C,
\end{equation}
for any $x$ as a member of $\mathbb{Q}^{n}_{p}$ and  all $\gamma$ in $\mathbb{Z}$. When there is a constant C which is positive we say that u lies in $W^{\infty}(\mathbb{Q}^{n}_{p})$, for which
\begin{equation}\label{C2e19}
|u(x)-u(y)|\leq C\frac{1}{log_{p}\left(p+min\{\|y\|_{p},\|x\|_{p}\}\right)},
\end{equation}
to any $x,y$ belongs to $\mathbb{Q}^{n}_{p}$.
Class $W^{\infty}_{0}(\mathbb{Q}^{n}_{p})$ is described as intersection of  $W_{0}(\mathbb{Q}^{n}_{p})$ and $W^{\infty}(\mathbb{Q}^{n}_{p})$.\\%$W^{\infty}_{0}(\mathbb{Q}^{n}_{p})=:W_{0}(\mathbb{Q}^{n}_{p})\cap W^{\infty}(\mathbb{Q}^{n}_{p})$

M is the Hardy-Littlewood maximal operator for a locally integrable function $f$ on $\mathbb{Q}^{n}_{p}$ in the following way:
\begin{equation*}
Mf(x)=\sup_{\gamma\in\mathbb{Z}}\frac{1}{p^{n\gamma}}\int_{B_{\gamma}}|f(\mathfrak{y})|d\mathfrak{y}.
\end{equation*}
Set $\mathcal{B}(\mathbb{Q}^{n}_{p})$ is of the form $u(\cdot)\in\mathfrak{G}(\mathbb{Q}^{n}_{p})$ meeting the boundedness condition for M between $L^{u(\cdot)}(\mathbb{Q}^{n}_{p})$ and $L^{v(\cdot)}(\mathbb{Q}^{n}_{p})$ (where $v$ is the Sobolev limiting exponent exponent, see (\ref{C5e1})).It is generally known that the analysis is significantly influenced by the Hardy-Littlewood maximal operator's boundedness on Lebesgue spaces.

Now, we will define p-adic variable exponent function spaces.

\begin{definition}\label{C1D4}
A function $f\in L^{u({\cdot})}_{loc}(\mathbb{Q}^{n}_{p})$ for $u({\cdot})\in\mathfrak{G}(\mathbb{Q}^{n}_{p})$ is  in $p$-adic $CMO^{u({\cdot})}(\mathbb{Q}^{n}_{p})$ with variable exponent if
\begin{equation*}
\|f\|_{CMO^{u({\cdot})}(\mathbb{Q}^{n}_{p})}=:\sup\limits_{\gamma\in\mathbb{Z}}\|\chi_{B_{\gamma}}\|^{-1}_{L^{u({\cdot})}(\mathbb{Q}^{n}_{p})}\|(f-f_{B_{\gamma}})\|_{L^{u({\cdot})}(\mathbb{Q}^{n}_{p})}<\infty,
\end{equation*}
where
\begin{equation*}
f_{B_{\gamma}}=\frac{1}{|B_{\gamma}|}\int_{B_{\gamma}}f(x)dx.
\end{equation*}
\end{definition}
If $u(x)=u$ is a constant, then $CMO^{u({\cdot})}(\mathbb{Q}^{n}_{p})$ equals $CMO^{u}(\mathbb{Q}^{n}_{p})$. We write $C^{u({\cdot})}=:CMO^{u({\cdot})}(\mathbb{Q}^{n}_{p})$ simply here and in the following.

\begin{definition}\label{C1D7}Let $\beta\in\mathbb{R}$, $0<m<\infty$, and $u(\cdot)\in\mathfrak{G}(\mathbb{Q}^{n}_{p})$. $\dot{K}^{\beta,m}_{u(\cdot)}(\mathbb{Q}^{n}_{p})$ is the homogeneous version of $p$-adic Herz space and its norm is given by
\begin{equation*}
	\dot{K}^{\beta,m}_{u(\cdot)}(\mathbb{Q}^{n}_{p})= \left \{ {g}\in L^{u(\cdot)} _{\mathrm{loc}}(\mathbb{Q}^{n}_{p}) : \| {g} \| _ {	 \dot{K}^{\beta,m}_{u(\cdot)}(\mathbb{Q}^{n}_{p})}  < \infty \right\},
	\end{equation*}
	where
	$$ \| {g} \| _{	\dot{K}^{\beta, m}_{u(\cdot)}(\mathbb{Q}^{n}_{p})} = \left ( \sum \limits _{\ell=-\infty} ^ {\infty} \| p^{\ell \beta} {g} \chi _\ell \| ^m
 _{L^{u(\cdot)}(\mathbb Q_p^n)}   \right )^ \frac{1}{m}.$$

\end{definition}

\begin{definition}\label{C1D8} Suppose $\beta\in\mathbb{R}$, $0<m<\infty$, $\lambda\in[0,\infty)$ and $u(\cdot)\in\mathfrak{G}(\mathbb{Q}^{n}_{p})$. $M\dot{K}^{\beta,\lambda}_{m,u(\cdot)}(\mathbb{Q}^{n}_{p})$ is the homogeneous version of $p$-adic Herz-Morrey space and its norm is given by
\begin{equation*}
	M\dot{K}^{\beta,\lambda}_{m,u(\cdot)}(\mathbb{Q}^{n}_{p})= \left \{ {g}\in L^{u(\cdot)} _{\mathrm{loc}}(\mathbb{Q}^{n}_{p}) : \| {g} \| _ {	 M\dot{K}^{\beta,\lambda}_{m,u(\cdot)}(\mathbb{Q}^{n}_{p})}  < \infty \right\},
	\end{equation*}
	where
	$$\| {g} \| _ {M\dot{K} ^{\beta, \lambda}_{m, u(\cdot)}(\mathbb{Q}^{n}_{p})} =\sup \limits_{k_0 \in \mathbb{Z}} 2^{-k_0 \lambda} \left( \sum \limits ^{k_0} _{\ell=-\infty}  \| p^{\ell\beta}{g} \chi _\ell  \| ^m _{L^{u(\cdot)}(\mathbb{Q}^{n}_{p})}\right)^\frac{1}{m}. $$

\end{definition}

\subsection{Fractional p-adic Hardy operator}

Hardy operators classical forms are defined by
\begin{equation*}
Hf(x)=:\frac{1}{x}\int^{x}_{0}f(\varsigma)d\varsigma,\hspace{0.4cm}H^{*}f(x)=:\int^{\infty}_{x}\frac{f(\varsigma)}{\varsigma}d\varsigma,\hspace{0.5cm}x>0,
\end{equation*}
for an integrable positive function $f$ on $\mathbb{R}^{+}$. Clearly, $\mathrm{H}$ and $\mathrm{H}^{*}$ satisfy the condition
\begin{equation*}
\int_{\mathbb{R}^{n}}g(\varsigma)Hf(\varsigma)d\varsigma=:\int_{\mathbb{R}^{n}}f(\varsigma)H^{*}g(\varsigma)d\varsigma.
\end{equation*}
For $1< u<\infty$, according to famous inequality for Hardy integrals \cite{160},
\begin{equation*}
\|Hf\|_{L^{u}(\mathbb{R}^{+})}\leq\frac{u}{u-1}\|f\|_{L^{u}(\mathbb{R}^{+})}.
\end{equation*}
Usefulness of Hardy integral inequalities in analysis and their applications have garnered considerable attention. With regards to their generalizations, variants, and applications there are numerous papers out there (cf. \cite{C5.5,C5.9,C5.10, C5.19,C2.25}).\\
On the other hand, Wu \cite{C5.21}defined fractional p-adic Hardy-type operator

\begin{definition}\label{D2}
Let $\mathrm{f},b\in L_{loc}(\mathbb{Q}^{n}_{p})$,\hspace{0.1cm}$0\leq\alpha<n.$ The fractional p-adic Hardy operators have the following definition:
\begin{equation}\label{C2e20}
H^{p}_{\alpha}f(x)=:\frac{1}{|x|^{n-\alpha}_{p}}\int_{|t|_{p}\leq|x|_{p}}f(t)dt,
\end{equation}
\begin{equation}\label{C2e21}
H^{p,*}_{\alpha}f(x)=:\int_{|t|_{p}>|x|_{p}}\frac{f(t)}{|t|^{n-\alpha}_{p}}dt,\hspace{0.2cm} x\in\mathbb{Q}^{n}_{p}\setminus\{0\},
\end{equation}
as well as their commutators
\begin{equation}\label{C2e22}
H^{p}_{\alpha,b}f=:bH^{p}_{\alpha}f-H^{p}_{\alpha}(bf),\hspace{0.2cm}H^{p,*}_{\alpha,b}f=:bH^{p,*}_{\alpha}f-H^{p,*}_{\alpha}(bf).
\end{equation}

\end{definition}
It is obvious that when $\alpha=0$, $H^{p}_{\alpha}$ turns into $H^{p}$.

Throughout this paper, many positive constants independent to primary variables will be denoted by the letter C. To keep things simple, we will use $\sum_{j=-\infty}^{\infty}f(x)\chi_{j}(x)=\sum_{j=-\infty}^{\infty}f_{j}(x)$.

Here are some lemmas regarding p-adic variable exponents that will help us to prove our main results.

\begin{lemma}\label{L1}\cite{C2.18}\hspace{0.1cm}Assume that $u\in\mathfrak{G}(\mathbb{Q}^{n}_{p})$ is an L-Lipschitz function for a value of $L\geq0,$ then $u\in W_{0}(\mathbb{Q}^{n}_{p}).$
\end{lemma}

\begin{lemma}\label{L2}\cite{C2.18}\hspace{0.1cm} Suppose $u\in W_{0}(\Omega^{n}_{p})$, where $\Omega^{n}_{p}\in\mathbb{Q}^{n}_{p}$ is a bounded set, then there arise an extension function $\tilde{u}\in W^{\infty}_{0}(\mathbb{Q}^{n}_{p})$ which is constant outside of some fixed ball.
\end{lemma}
\begin{lemma}\label{L5}\cite{C2.18}\hspace{0.1cm} Let $u(\cdot)\in W^{\infty}_{0}(\mathbb{Q}^{n}_{p})$. Then,
\begin{equation*}%\label{C2e30}
\|\chi_{B_{k}}\|_{L^{u(\cdot)}(\mathbb{Q}^{n}_{p})}\leq Cp^{kn/u(x,k)},
\end{equation*}
where
\begin{equation*}%\label{C2e31}
 u(x,k)=:
\begin{cases}
    u(x),&  k<0,\\
    u(\infty),              &  k\geq0.
\end{cases}
\end{equation*}
\end{lemma}

Lemma 2.2 in \cite{C2.29} is extended to the p-adic variable exponent central BMO space in the following Lemma.
\lemma\label{L3} \normalfont Let $g\in C^{u(\cdot)}$ and $m,l\in \mathbb Z,$ then \begin{equation}\label{C2eN1}|g(x)-g_{B_m}|\le|g(x)-g_{B_l}|+p^{n}|l-m|\|g\|_{C^{u(\cdot)}}.\end{equation}

\section{Variable $p$-adic Herz Space Estimates for Hardy Operators}
The findings of this section present the continuity characteristics about $H^{p}_{\alpha}$, $H^{p,*}_{\alpha}$, $H^{p}_{\alpha,b}$, and $H^{p,*}_{\alpha,}$, which are all associated with the variable exponent $p$-adic Herz space.
\begin{theorem}\label{C5T1}
Let $0<m_{1}\leq m_{2}<\infty$, $u(\cdot)\in\mathfrak{G}(\mathbb{Q}^{n}_{p})$, $0<\alpha<\min\{\frac{n}{u_{+}},\frac{n}{v^\prime_+}\}$ and $-\frac{n}{v_+}<\beta<\frac{n}{u_-^\prime}$. Defined $v(\cdot)$ by
\begin{equation}\label{C5e1}\frac{1}{v(\cdot)}=\frac{1}{u(\cdot)}-\frac{\alpha}{n},\end{equation} then both $H^{p}_{\alpha}$ and $H^{p,*}_{\alpha}$ map $\dot{K}^{\beta,m_{2}}_{v(\cdot)}(\mathbb{Q}^{n}_{p})$ into $\dot{K}^{\beta,m_{1}}_{u(\cdot)}(\mathbb{Q}^{n}_{p})$.
\end{theorem}
From the above theorem, if $\alpha=0$, then the following result is true.
\begin{corollary}\label{C5C1}
Let $0<m_{1}\leq m_{2}<\infty$, $u\in\mathfrak{G}(\mathbb{Q}^{n}_{p})$, and $-\frac{n}{u_+}<\beta<\frac{n}{u_-^\prime}$. Then both $H^{p}$ and $H^{p,*}$ map $\dot{K}^{\beta,m_{2}}_{u(\cdot)}(\mathbb{Q}^{n}_{p})$ into $\dot{K}^{\beta,m_{1}}_{u(\cdot)}(\mathbb{Q}^{n}_{p})$.
\end{corollary}

The next result gives the continuity of commutators of $p$-adic Hardy-type operators on $p$-adic variables exponent Herz space.
\begin{theorem}\label{C5T2}
Let $0<m_{1}\leq m_{2}<\infty$, $b\in C^{u^\prime(\cdot)}\cap C^{v(\cdot)}$, $v(\cdot)\in W^{\infty}_{0}(\mathbb{Q}^{n}_{p})$, $0<\alpha<\min\{\frac{n}{u_{+}},\frac{n}{v^\prime_+}\}$ and $-\frac{n}{v_+}<\beta<\frac{n}{u_-^\prime}$. Defined $v(\cdot)$ by
\begin{equation}\label{C5e1}\frac{1}{v(\cdot)}=\frac{1}{u(\cdot)}-\frac{\alpha}{n},\end{equation} then both $H^{p}_{\alpha,b}$ and $H^{p,*}_{\alpha,b}$ map $\dot{K}^{\beta,m_{2}}_{v(\cdot)}(\mathbb{Q}^{n}_{p})$ into $\dot{K}^{\beta,m_{1}}_{u(\cdot)}(\mathbb{Q}^{n}_{p})$.
\end{theorem}
The following corollary holds if $\alpha=0$ in the preceding theorem.
\begin{corollary}\label{C5C2}
Let $0<m_{1}\leq m_{2}<\infty$, $b\in C^{u^\prime(\cdot)}\cap C^{u(\cdot)}$, $u\in\mathfrak{G}(\mathbb{Q}^{n}_{p})$, and $-\frac{n}{u_+}<\beta<\frac{n}{u_-^\prime}$. Then both $H^{p}_{b}$ and $H^{p,*}_{b}$ map $\dot{K}^{\beta,m_{2}}_{u(\cdot)}(\mathbb{Q}^{n}_{p})$ into $\dot{K}^{\beta,m_{1}}_{u^\prime(\cdot)}(\mathbb{Q}^{n}_{p})$.
\end{corollary}

\section{Variable Morrey-Herz Estimates for $p$-adic Hardy-type Operators and Commutators}
This section proves the boundedness of $H^{p}_{\alpha}$, $H^{p,*}_{\alpha}$, $H^{p}_{\alpha,b}$ and $H^{p,*}_{\alpha,b}$ on Morrey-Herz type spaces. Here $f_{i}=f(\chi_{i})$ remains the same as used in previous section for any $i\in\mathbb{Z}$.
\begin{theorem}\label{C5T3}
Let $0<m_{1}\leq m_{2}<\infty$, $u(\cdot)\in \mathfrak{G}(\mathbb{Q}^{n}_{p})$, $0<\alpha<\min\{\frac{n}{u_{+}},\frac{n}{v^\prime_+}\}$ and $\lambda-\frac{n}{v_-}<\beta<\frac{n}{u^\prime_-}+\lambda$. Defined $v(\cdot)$ by
\begin{equation}\label{C5e1}\frac{1}{v(\cdot)}=\frac{1}{u(\cdot)}-\frac{\alpha}{n},\end{equation} then both $H^{p}_{\alpha}$ and $H^{p,*}_{\alpha}$ map $M\dot{K}^{\beta,\lambda}_{m_{2},v(\cdot)}(\mathbb{Q}^{n}_{p})$ into $M\dot{K}^{\beta,\lambda}_{m_{1},u(\cdot)}(\mathbb{Q}^{n}_{p})$.
\end{theorem}
If $\alpha=0$, then the following is true:
\begin{corollary}\label{C5C3}
Let $0<m_{1}\leq m_{2}<\infty$, $u\in\mathfrak{G}(\mathbb{Q}^{n}_{p})$, and $-\frac{n}{u_+}<\beta<\frac{n}{u^\prime_-}$. Then both $H^{p}$ and $H^{p,*}$ map $M\dot{K}^{\beta,\lambda}_{m_{2},u(\cdot)}(\mathbb{Q}^{n}_{p})$ into $M\dot{K}^{\beta,\lambda}_{m_{1},u(\cdot)}(\mathbb{Q}^{n}_{p})$.
\end{corollary}

\begin{theorem}\label{C5T4}
Let $0<m_{1}\leq m_{2}<\infty$, $b\in C^{u^\prime(\cdot)}\cap C^{v(\cdot)}$, $u(\cdot)\in \mathfrak{G}(\mathbb{Q}^{n}_{p})$, $0<\alpha<\min\{\frac{n}{u_{+}},\frac{n}{v^\prime_+}\}$ and $\lambda-\frac{n}{v_-}<\beta<\frac{n}{u^\prime_-}+\lambda$. Defined $v(\cdot)$ by
\begin{equation}\label{C5e1}\frac{1}{v(\cdot)}=\frac{1}{u(\cdot)}-\frac{\alpha}{n},\end{equation} then both $H^{p}_{\alpha,b}$ and $H^{p,*}_{\alpha,b}$ map $M\dot{K}^{\beta,\lambda}_{m_{2},v(\cdot)}(\mathbb{Q}^{n}_{p})$ into $M\dot{K}^{\beta,\lambda}_{m_{1},u(\cdot)}(\mathbb{Q}^{n}_{p})$.
\end{theorem}
The logical consequence of $\alpha=0$ is as follows:
\begin{corollary}\label{C5C4}
Let $0<m_{1}\leq m_{2}<\infty$, $b\in C^{u^\prime(\cdot)}\cap C^{u(\cdot)}$, $u(\cdot)\in\mathfrak{G}(\mathbb{Q}^{n}_{p})$, and $\lambda-\frac{n}{u_-}<\beta<\frac{n}{u^\prime_-}+\lambda$. Then both $H^{p}_{b}$ and $H^{p,*}_{b}$ map $M\dot{K}^{\beta,\lambda}_{m_{2},u(\cdot)}(\mathbb{Q}^{n}_{p})$ into $M\dot{K}^{\beta,\lambda}_{m_{1},u^\prime(\cdot)}(\mathbb{Q}^{n}_{p})$.
\end{corollary}

\end{document}